\documentclass[a4paper, 12pt]{article}

\usepackage{amssymb,amsbsy,amsmath,amsfonts,amssymb,amscd}
\usepackage{latexsym,eucal,exscale,epic,eepic,epsfig,xcolor}

\newtheorem{theorem}{Theorem}[section]
\newtheorem{lemma}[theorem]{Lemma}
\newtheorem{proposition}[theorem]{Proposition}
\newtheorem{corollary}[theorem]{Corollary}

\newtheorem{conjecture}[theorem]{Conjecture}

\newcommand{\ds}{\displaystyle}

\newcommand{\ofr}{{\mathfrak o}}
\newcommand{\pfr}{{\mathfrak p}}

\newcommand{\lra}{\longrightarrow}
\newcommand{\noi}{\noindent}

\newcommand{\St}{{\mathbf S}{\mathbf t}}

\newcommand{\ZZ}{\mathbb Z}

\newcommand{\Sm}{\mathcal S}
\newcommand{\CC}{\mathbb C}

\newcommand{\HH}{{\mathbb H}}
\newcommand{\HIw}{{\mathcal H}_{\rm I}}
\newcommand{\Iw}{{\rm I}}
\newcommand{\VV}{\mathcal V}
\newcommand{\FF}{\mathcal F}

\newcommand{\Ch}{\rm Ch}

\newcommand{\bM}{\mathbf M}

\newcommand{\WAff}{W^{\rm Aff}}
\newcommand{\SL}{{\rm SL}}

\title{On the distinction of Iwahori-spherical discrete series representations}
\author{Paul  Broussous}

\date{\today}

\begin{document}
\maketitle

\begin{abstract}
 Let $E/F$ be a quadratic unramified extension of non-archimedean local fields and $\mathbb H$ a simply connected semisimple algebraic group defined and split over $F$. We establish general results (multiplicities, test vectors) on $\HH (F)$-distinguished Iwahori-spherical representations of $\HH (E)$. For discrete series Iwahori-spherical representations of $\HH (E)$, we prove a numerical criterion of $\HH (F)$-distinction. As an application, we classify the $\HH (F)$-distinguished discrete series representations of $\HH (E)$ corresponding to degree $1$ characters of the Iwahori-Hecke algebra. 
\end{abstract}

 \tableofcontents

 \bigskip

\section{Introduction} 

 This work\footnote{AMS Classification: 22E50, 22E35, 11F70, 20C08. \\ Keywords: reductive groups over local fields, relative Langlands program, reductive symmetric spaces, Iwahori-Hecke algebras.} is a contribution to the study of distinguished complex representations for Galois symmetric spaces over  a non-archimedean local field. It may be viewed as a continuation of \cite{BC} where the authors prove a particular case of Prasad's conjecture on the distinction of the Steinberg representation (\cite{Pr2}, Conjecture 3, p. 77). However even though the Bruhat-Tits building will be used here as a geometrical tool as in \cite{BC}, our fundamental techniques will be that of integration of matrix coefficients, in the spirit of \cite{Br}.

  We restrict ourselves to the easy case where the symmetric space is of the form ${\mathbb H}(E)/{\mathbb H}(F)$, where $E/F$ is a quadratic unramified extension of non-archimedean local fields and $\mathbb H$ a connected  simple group defined and split over $F$.  For simplicity sake we shall assume that $\mathbb H$ is simply connected. Indeed in that case the set of Iwahori subgroups of ${\mathbb H}(E)$ identifies with the set of chambers in the Bruhat-Tits building of ${\mathbb H}(E)$, and this  makes our calculations less technical.  However our arguments should adapt easily to the non-simply connected case.

 Recall that a smooth complex representation $(\pi ,\VV )$ of $\HH (E)$ is said to be $\HH (F)$-distinguished if the intertwining space ${\rm Hom}_{\HH (F)} \, (\VV ,\CC )$ is non-zero, where $\CC$ is a complex line with the trivial action of $\HH (F)$. In this article we study the distinction of {\it Iwahori-spherical} irreducible representations, i.e. representations $(\pi ,\VV )$ such that $\VV^I\not= 0$, where $I$ is some fixed Iwahori subgroup of $\HH (E)$. If $(\pi ,\VV )$ is such a  representation, the Iwahori-Hecke algebra $\HIw$ of $\HH (E)$ relative to $I$ acts on $M=\VV^I$ via an irreducible representation $(r,M)$. 
\smallskip

 Our first result is the following.

\begin{theorem} (Theorem  \ref{Multiplicite}) Let $(\pi ,\VV )$ be an Iwahori-spherical irreducible representation of $\HH (E)$. 
 Let $J$ be any ${\rm Gal}(E/F)$-stable Iwahori subgroup of $\HH(E)$. Then the obvious restriction map ${\rm Hom}_{\HH (F)}\, (\VV ,\CC )\lra {\rm Hom}_{\CC}\, (\VV^J ,\CC )$ is injective. 

As a consequence, If $\Lambda \in {\rm Hom}_{\HH (F)} (\VV ,\CC )$ is non-zero then for any  ${\rm Gal}(E/F)$-stable Iwahori subgroup $J$ of $\HH (E)$, there exists  a vector $v\in \VV^J$ such that $\Lambda (v)\not= 0$. 
\end{theorem}

\noi The proofs of these results are based on a transitivity property of the action of $\HH (F)$ on the set of chambers of the Bruhat-Tits building of $\HH (E)$, which was a crucial ingredient of \cite{BC}. 
 \medskip

 We next assume moreover that $(\pi ,\VV )$ is a discrete series representation. Using the fact due to Zhang \cite{Zh} that $\HH (E )/\HH (F)$ is a very strongly discrete symmetric space, we establish a second result. 

\begin{theorem} (Theorem \ref{Bitest}) Let $(\pi ,\VV )$ be an irreducible Iwahori-spherical discrete series representation of $\HH (E)$. Let $J$ be some ${\rm Gal}(E/F)$-stable Iwahori subgroup of $\HH (E)$ and $\nu$ denote a Haar measure on $\HH (F)$. 

Then $(\pi ,\VV )$ is $\HH (F)$-distinguished if, and only if, there exist $v\in \VV^J$, ${\tilde
 v}\in {\tilde \VV}^J$, such that 
$$
\int_{\HH (F)} c_{v,{\tilde v}}(h)\, d\nu (h)\not= 0\ .
$$
\noi where $\tilde\VV$ is the contragredient representation and $c_{v, {\tilde v}}$ is the coefficient of $\pi$ attached to $v$ and $\tilde v$. 
\end{theorem}

To turn the previous theorem into an effective criterion, we use the notion of {\it generalized Poincar\'e series} attached to an Iwahori-Hecke algebra introduced by Gyoja \cite{Gy}. Let $(W,S)$ be a Coxeter system with associated Hecke algebra ${\mathcal H}$. Let $S_i$, $i=1,...,m$ be the intersections of the conjugacy classes of $W$ with $S$. For each representation $(r,M)$ of ${\mathcal H}$, Gyoja attaches a power series $L(t_1 ,...,t_m ,r)\in {\rm End}(M)\otimes_\CC \CC [[t_1 ,...,t_m ]]$.  Now assume that ${\mathcal H}=\HIw$ is the Iwahori-Hecke algebra of $\HH (E)$ relative to $I$ and $(r,M)$ is the representation associated to the $\HIw$-module $\VV^I$, for a representation $(\pi , \VV )$ of $\HH (E)$.

\begin{theorem} \label{criterion}  (Theorem \ref{Criterion}) Let $(\pi ,\VV )$ be an irreducible Iwahori-spherical discrete series representation of $\HH (E)$. Then $\pi$ is $\HH (F)$-distinguished if, and only if $L(1/q_o ,...,1/q_o ,r)\in {\rm End} (M)$ is non-zero, where $q_o$ is the size of the residue field of $F$. 
\end{theorem}

 Our last result is concerned with the case where the discrete series $(\pi ,\VV )$ corresponds to a degree $1$ character of $\HIw$. Those representations {\color{blue}} were classified by A. Borel in \cite{Bo}. It is by using techniques borrowed from {\it loc. cit} that we prove the following:  the value of $L (t_1 ,...,t_m ,r)$ occuring in Theorem \ref{criterion} can be calculated using a Macdonald's Poincar\'e series  (cf. \cite{Mc}).  Using Macdonald's formulas, we obtain:

\begin{theorem} (Theorem \ref{Main}) Let $(\pi ,\VV )$ be an irreducible Iwahori-spherical discrete series representation of $\HH (E)$ such that $\VV^I$ is $1$-dimensional. If $\HH$ is not of type $G_2$, then $\pi$ is $\HH (F)$-distinguished if, and only if, it is the Steinberg representation. If $\HH$ has type $G_2$ there are two representations satisfying the assumptions of the theorem: the Steinberg representation and another one; there are both $\HH (F)$-distinguished. 
\end{theorem}

\noi In other words, the Steinberg representation of $\HH (E)$ is the unique representation satisfying the assumption of the theorem and which is $\HH (F)$-distinguished, except in type $G_2$ where there is another one. 
\medskip

 \noi The case of the Steinberg representation (originally D. Prasad's conjecture) was already established in \cite{BC} (where $\HH$ is not necessarilly assumed simply connected) and in characteristic $0$ Prasad's conjecture was entirely proved by Beuzart-Plessis \cite{Be}. 
\medskip

 Macdonald's formulas \cite{Mc} allow only to calculate Gyoja $L$-functions $L(t_1 ,...,t_m ,r)$ when $r$ is a degree $1$ character, this is why we are not able to give a classification of distinguished discrete series representations when ${\rm dim}\, \VV^I >1$. 
\medskip

In \cite{Pr}, Dipendra Prasad states a conjecture predicting the distinction of irreducible smooth complex representations in the case of a general Galois symmetric space $\HH (E)/\HH (F)$  (Conjecture 2 of {\it loc. cit.}).  An irreducible smooth representation $\pi$ of $\HH (E)$ conjecturally has a Langlands parameter $\phi_\pi$~: $W_E \times \SL_2 (\CC )\lra {}^L G$, where the Langlands dual ${}^L G$ is a semidirect product $W_E\ltimes \widehat{G}$, $\widehat{G}$ being the complex dual of $G=\HH (E)$.  Roughly speaking, the point of Prasad's conjecture is to relate the $\HH (F)$-distinction of  $\pi$ to properties of the parameter $\varphi_\pi$. A precise formulation of this conjecture in our case may be stated as follows (this formulation as well as the two following remarks are due to Prasad himself).

\begin{conjecture} (Prasad) Let $E/F$ be a quadratic unramified extension of nonarchimedean local fields. Let $\HH$ be a simply connected split simple group
  over $F$ if $-1$ belongs to the Weyl group of $\HH$, and the unique quasi-split (but not split) group
  over $F$ which splits over $E$ otherwise. (This ensures that the opposition group $H^{\rm op}$ of $H=\HH (F)$, defined in \cite{Pr}, 
is split over $F$.) A Iwahori-spherical discrete series representation  $\pi$ of $G=\HH (E)$   has, by the work of Kazhdan-Lusztig, cf. \cite{KL}, an $\rm L$-parameter  
$\phi: W_E \times \SL_2(\CC ) \rightarrow {}^LG$ which on $W_E$ factors through an unramified extension of $E$.
Assume that the group of connected components of the centralizer of the corresponding unipotent element
in $\widehat{G}$ is  an elementary abelian $2$-group (which is known to be the case for all
groups of type $B_n,C_n,D_n$). Then $\pi$
is distinguished by $H= \HH (F)$ if and only if

\begin{enumerate}
\item The restriction of $\phi$ to $W_E$ is trivial,
\item The representation $\pi$ is generic. (Under condition (1), genericity should not depend
  on the choice of the Whittaker datum.)
\end{enumerate}

Further, if the conditions (1) and (2) are satisfied, then ${\rm dim}  \, {\rm Hom}_{H}\,  (\pi, \CC )$ is the  order of the centralizer of $\phi(\SL_2(\CC ))$ in $ \hat{G}$,
which by assumption  is a finite elementary abelian 2-group.

\end{conjecture}

\noi {\bf Remark 1.} We relate the above Conjecture to that in \cite{Pr}.  In the present case, $W_F \times \SL_2(\CC )$ and
  $W_E\times \SL_2(\CC)$ can be replaced by $\langle {\rm Fr}^{\ZZ} \rangle \times \SL_2(\CC )$ and
  $\langle {\rm Fr}^{2\ZZ} \rangle\times \SL_2(\CC)$, where ${\rm Fr}$ is a Frobenius element of $W_F$.
  Because of the hypothesis in the Conjecture, $\pi_0 (Z_{\widehat{G}}(\phi ( \SL_2(\CC ))))$ being an elementary abelian $2$-group,   $\langle {\rm Fr}^{2\ZZ} \rangle$ must land inside the connected component of identity of
$Z_{\widehat{G}}(\phi (\SL_2(\CC )))$, hence by the discreteness
  of the parameter $\phi$,  $Z_{\widehat{G}}(\phi (\SL_2(\CC )))$
  must be finite and therefore an elementary abelian $2$-group, and the image of  $\langle {\rm Fr}^{2\ZZ} \rangle$ under $\phi$ must be trivial  which is the conclusion in part
  (1) of the Conjecture. Further, the natural map from the component group
  of a lifted parameter $\tilde{\phi}$ for $H^{\rm op}$ to the component group for $\phi$ is clearly an
  isomorphism, hence by the Conjecture of \cite{Pr},
  the only $\pi$ corresponding to the trivial character of the component group can be
  distinguished which is the conclusion in part (2) of the Conjecture.
\medskip

\noi {\bf Remark 2.} Here is the example of $G_2$ which appears in this work but for which the hypothesis in the Conjecture above does not hold. There is a discrete series representation $\pi$ of $G_2 (E)$ which has a one dimensional fixed vector
  for an Iwahori subgroup, which by Reeder \cite{Re}, page 480, corresponds to a subregular unipotent element with component group isomorphic to $S_3$ with the representation
  of the component group, the $2$-dimensional irreducible representation $\rho$ of $S_3$. It is thus a non-generic representation of $G_2$ which is distinguished
  by $G_2(F)$ (by Theorem \ref{Main} of this paper). In this case, the representation $\rho$ of $S_3$ has a fixed vector under any transposition in $S_3$, thus confirming the general conjecture
  of \cite{Pr}  in an example with non-abelian component group.

 It would be nice to check Prasad's conjecture in our case. In this aim, one would have to determine the  Galois parameters of our representations, which are  not known in general.

\medskip

The article is organized as follows. After introducing the general notation in {\S}1, in {\S}2 we recall the models of Iwahori-spherical representations introduced in \cite{Bo} and explain how to deduce a very simple formula for the Iwahori-spherical coefficients. The bound on the dimension of the space of $\HH (F)$-invariant linear forms is given in {\S}3 and the fact that test-vectors can be chosen Iwahori-spherical is proved in {\S}4. In {\S}5 we recall the definition of Gyoja's generalized Poincar\'e series and we give Macdonald's formulas. Finally {\S}6 is devoted to establishing the numerical criterion in terms of Gyoja's Poincar\'e series values and to applying it to the case of degree $1$ characters of the Iwahori-Hecke algebra. 
\medskip

I warmly thank Nadir Matringe and Dipendra Prasad for their help in writing this article.

\section{Notation} Throughout the article we use the following notation. 

\medskip

 If $X$ is a set, we denote by $\vert X\vert$ its cardinal. If a group $\Gamma$ acts on $X$ and $\Gamma '$ is a subset of $\Gamma$, $X^{\Gamma '}$ denotes the subset of those $x\in X$ that are fixed by  $\Gamma '$.  
\smallskip

 We let $F$ denote a non-archimedean, non-discrete,  locally compact  field. We denote by $\ofr_F$ its ring of integers, by $\pfr_F$ the maximal ideal of $\ofr_F$, and we set $k_F =\ofr_F /\pfr_F$ (residue field), $q_0 =\vert k_F\vert  =p^r$, $p$ being the characteristic of $k_F$ (residue characteristic). 
\smallskip

 We fix an {\it unramified quadratic extension} $E/F$ and use the obvious pieces of notation $\ofr_E$, $\pfr_E$, $k_E$. We denote by $\theta$ the generator of the Galois group ${\rm Gal}(E/F)$. We set $q =\vert k_E\vert$, so that $q=q_0^2$. 
\smallskip

 We fix a connected  simple $F$-algebraic group $\mathbb H$ assumed to be {\it split over $F$ and simply connected} and we write $H={\mathbb H}(F)$ for the locally compact group of its $F$-rationnal points. We denote by ${\mathbb G}={\rm Res}_{E/F}{\mathbb H}$ the $F$-algebraic group obtained by Weil's restriction of scalars. It is  semisimple and simply connected. We set $G={\mathbb G}(F)={\mathbb H}(E)$. 

 The Galois group ${\rm Gal}(E/F)$ acts on $\mathbb G$ by $F$-rational automorphisms of algebraic group. We still denote by $\theta$ the action of $\theta$ on $\mathbb G$ and $G$ so that ${\mathbb H} = {\mathbb G}^\theta$ and $H=G^\theta$. The pair $(G,H)$ (or the $G$-set $G/H$) is called a {\it (Galois) reductive symmetric space}.
\smallskip

We let $X_F$ (resp. $X_E$) denote the Bruhat-Tits building of $\mathbb H$ over $F$ (resp. over $E$). Recall that $X_E$ identifies canonically with the Bruhat-Tits building of ${\rm Res}_{E/F}({\mathbb H})$ over $F$. As simplicial complexes $X_F$ and $X_E$ have dimension $d$, the $F$-rank of $\mathbb H$. As a particular case of building functorialities, there is an injection $j$ : $X_F\lra X_E$ satisfying the following properties:
\smallskip

(a) $j$ is ${\rm Gal}(E/F)\ltimes H$-equivariant, 

(b) the image of $j$ is $X_E^{{\rm Gal}(E/F)}$, 

(c) $j$ is simplicial and maps chambers to chambers. 
\smallskip

 We identify $X_F$ as a subsimplicial complex of $X_E$ viewing $j$ as an inclusion. 
\smallskip

We fix a maximal $F$-split torus $\mathbb T$ of $\mathbb H$ and set $T_F = {\mathbb T}(F)$, $T_E ={\mathbb T}(E)$.  We denote by $\Phi =\Phi ({\mathbb T}, {\mathbb H})$ the root system of $\mathbb T$ in $\mathbb H$. Let $W^{\rm Aff}_F = N_{H} ({T} )/T_F^o$ be the affine Weyl group of $H$ relative to $\mathbb T$, where $N_{ H} ( T_F )$  is the normalizer of $\mathbb T$ in $H$ and $T_F^o$ the maximal compact subgroup of $T_F$. The group $W^{\rm Aff}_E$ is defined in the same way and the inclusion $N_H (T_F )\subset N_G (T_E )$ induces an isomorphism of groups $W_F^{\rm Aff} \simeq W_E^{\rm Aff}$.  In the sequel we abbreviate $W^{\rm Aff} =W_F^{\rm Aff}=W_E^{\rm Aff}$. 

 Let ${\mathcal A}_F$ (resp. ${\mathcal A}_E$ ) denote the appartment of $X_F$ (resp. $X_E$) associated to $\mathbb T$. It follows from the construction of $j$ that ${\mathcal A}_F ={\mathcal A}_E$. Fix a chamber $C_o$ of ${\mathcal A}_F$ and let $I^o_F$ (resp. $I_E^o$)  denote the Iwahori subgroup of $H$ (resp. the ${\rm Gal}(E/F)$-stable Iwahori subgroup  of $G$) fixing $C_o$. Since $\mathbb H$ is simply connected these Iwahori subgroups are the (global) fixators of $C_o$.  Let $S\subset W_F^{\rm Aff}$ be the set of involutions corresponding to reflections in ${\mathcal A}_F$ relative to the walls of $C_o$. Then, since $\HH$ is simply connected,  $(I_F^o ,N_H (T_F ))$ (resp. $(I_E^o ,N_G (T_E ))$) is a Tits system (or BN-pair) in $H$ (resp. in $G$) with Coxeter system $(W^{\rm Aff},S)$. We denote by $l$ the length function of $W$ relative to the generative subset $S$. Recall that 
$$
\vert I_E^o w I_E^o /I_E^o \vert =q^{l(w)}, \ w\in \WAff\ .
$$

Recall that the vertices of $X_E$ are naturally labelled by the elements of $S$. If $C$ is a chamber of $X_E$, the  wall of $C$ which does not contain the unique vertex of type $s$ is said to be of type $s$. If $C$ and $D$ are two  chambers of $X_E$, we write $D\sim_s C$ if $D$ and $C$ are adjacent and  if the wall $D\cap C$ is of type $s$.

\smallskip

 We let $\HIw$ denote the Iwahori-Hecke algebra of $G$ relative to our choice $I_E^o$ of Iwahori subgroup. Recall that as a complex vector space $\HIw$ is the set of complex functions of $G$ which are bi-invariant by $I_E^o$ and have compact support. It is endowed with the convolution product defined by 
$$
u\star v (g)=\int_G u(gx^{-1})v(x)\, d\mu (x)=\sum_{x\in I_E^o \backslash G} u(gx^{-1})v(x), \ g\in G, \ u,v\in \HIw
$$
\noi where $\mu$ is the Haar measure on $G$ normalized by $\mu (I_E^o )=1$.  Recall that $(e_w )_{w\in \WAff}$, is a $\CC$-basis of $\HIw$, where for $w\in \WAff$, $e_w$ denote the caracteristic function of $I_E^o wI_E^o$. Recall that we have the classical relations:
\begin{align*}
(e_s +1)\star (e_s -q)=0, \ s\in S\\
e_w \star e_{w'} = e_{ww'}\ {\rm if}\ l(ww')=l(w)+l(w'), \ w, \ w'\in W^{\rm Aff}
\end{align*}
\noi and that the $e_s$, $s\in S$, generate $\HIw$ as an algebra. 

 If $f\in \HIw$, we define $\check{f}$ by $\check{f} (x)=f( x^{-1} )$, $x\in G$. The map $f\mapsto \check{f}$ is an automorphism of the vector space $\HIw$ satisfying $ (f\star g)\check{} =\check{g} \star \check{f}$, $f,g\in \HIw$. If $(r,M)$ is a complex representation of $\HIw$, we define its contragredient $({\tilde r}, {\tilde M})$ in the dual ${\tilde M}$ of $M$ by ${\tilde r}(f)={}^t r(\check{f} )$, $f\in \HIw$, where ${}^t$ denotes the transposition.  
 \smallskip

In this article, representations of $G$, $H$ or $\HIw$ are always assumed to be in complex vector spaces. 
If $(\pi ,\VV )$ is a smooth representation of $G$, $({\tilde \pi}, {\tilde \VV })$ denotes its contragredient (i.e. its smooth dual). If $W$ is a complex vector space, $W^* ={\rm Hom}_\CC (W,\CC )$ denotes its algebraic dual.

\section{Iwahori-spherical representations: models and coefficients}

In this section we abbreviate $I=I_E^o$. 
\smallskip

Let $\Sm (G)$ denote the abelian category of smooth representations of $G$ and $\Sm (G)_{\rm Iw}$ denote the full subcategory of those representations $(\pi , \VV )$ that are generated by the fixed vector set $\VV^I$. 
Let $\HIw{\rm -Mod}$ denote the category of left $\HIw$-modules. 
If $(\pi ,\VV )$ is an object of $\Sm (G)_\Iw$, then $\VV^I$ is naturally a left $\HIw$-module, so that we have a well defined functor 
$$
{\mathbf M} \ : \ \Sm (G)_\Iw\lra \HIw{-}{\rm Mod}, \ (\pi ,\VV )\mapsto \VV^I\ .
$$
\noi It is now a well known fact (\cite{Bo},\cite{Ca}) that $\mathbf M$ is an equivalence of categories; in other words the trivial character of $I$ is a {\it type} for $G$ in the sense of \cite{BK2}. Moreover if $(\pi ,\VV )$ is admissible, we have $\bM ({\tilde \VV} )= \bM (\VV )\, \tilde{}$ : the functor commutes with the operation of taking the contragredient.

 A pseudo-inverse $\mathbf V$ for $\mathbf M$ may be constructed as follows (cf. \cite{Bo}{\S}4). Let us fix a left $\HIw$-module $M$ viewed as a representation $(r,M)$ of $\HIw$. Let $C_c (G/I)$ be the space of  functions on $G$ which are right $I$-invariant and have compact support. This is naturally a right $\HIw$-module and a smooth representation of $G$ (by left translation) and both structures are compatible. Hence we may form the tensor product ${\mathbf V}(M) := C_c (G/I)\otimes_{\HIw} M$.  Of course we have an isomorphism of left $\HIw$-modules $M\simeq {\mathbf V}(M)^I$; it is explicitly given by  
$m\mapsto e_1 \otimes m$.

\medskip

We shall need two other functorial constructions  $\mathbf P$, ${\mathbf P}^o$,   considered by Borel in loc. cit. {\S}2.  If $(r,M)$ is a representation of $\HIw$ we set 
$$
{\mathbf P}^o (M)=\{ f\in C(G/I,M)\ ; \ f\star r(u) =r(u).f , \ u\in \HIw\}\ ,
$$
\noi where $C(G/I ,M)$ denotes the space of $M$-valued functions on $G$ which are right $I$-invariant. This space is acted upon by $G$ via left translation, but the obtained representation of $G$ is not smooth in general; we denote by ${\mathbf P}(M)\subset {\mathbf P}^o (M)$ the subspace of smooth vectors. By \cite{Bo} Proposition 2.4, the map 
$$
\nu_o \ : \ {\mathbf P} (M)^I ={\mathbf P}^o (M)^I\lra M\ , \ f\mapsto f(1)\ ,
$$
\noi is an isomorphism of vector spaces. Moreover, for $m\in M$,  the element  $f_m\in {\mathbf P}(M)^I$ mapped onto $m$ by $\nu_o$ is the function $f_e$ on $G$ mapping any element of $IwI$, $w\in\WAff$, to $q_w^{-1}\, r(\check{e_w}).m$, symbolically written by Borel as
\begin{equation}\label{FormeIwahori}
f_m = \sum_{w\in \WAff} e_w\ q_w^{-1}\, r(\check{e_w}).m\ .
\end{equation}
 
The tensor product $C_c (G/I)\otimes_{\HIw} M$ may be viewed as the quotient of $C_c (G/I)\otimes_\CC M$ by the $G$-invariant subspace generated by $(f\star u)\otimes m-f\otimes r(u).m$, $f\in C_c (G/I)$, $u\in \HIw$, $m\in M$. We may identify $C_c (G/I)\otimes_\CC M$ with the space $C_c (G/I,M)$ of $M$-valued functions on $G$ which are right $I$-invariant. The bilinear $G$-invariant map
\begin{equation}\label{Pairing}
{\mathbf P}^o ({\tilde M})\times C_c (G/I,M)\lra \CC, \ (f,g)\mapsto \sum_{x\in G/I} \langle f(x), g(x)\rangle_M
\end{equation}
induces an well-defined $G$-invariant bilinear map
$$
{\mathbf P}^o ({\tilde M})\times {\mathbf V}(M)\lra \CC
$$
whence a $G$-intertwining operator
$$
\psi_M\ :\ {\mathbf P}^o ({\tilde M})\lra {\mathbf V}(M)^*
$$
\noi where $*$ denotes an algebraic dual.  From \cite{Bo} Proposition 2.6, we have the following.

\begin{proposition} \label{Dual}Let $(r,M)$ be a representation of $\HIw$. Then the intertwining operator $\psi_M\ :\ {\mathbf P}^o ({\tilde M})\lra {\mathbf V}(M)^*$ is an isomorphism and restricts to an isomorphism
${\mathbf P} ({\tilde M})\simeq {\mathbf V}(M)\tilde{}$. 
\end{proposition}

 Observe that if $M$ is finite dimensional, ${\mathbf V}(M)$ is admissible and we have ${\mathbf V}(M)\tilde{} \simeq {\mathbf V}({\tilde M})$. So exchanging the roles of $M$ and $\tilde M$, we obtain ${\mathbf P}(M)\simeq {\mathbf V}(M)$ as $G$-modules.
\medskip

 Let $(r,M)$ be a finite dimensional representation of $\HIw$ and set $(\pi ,\VV )={\mathbf V}(r,M)$. We make the canonical identifications $\VV^I =M$ and ${\tilde \VV}^I ={\tilde M}$. Recall that for $m\in M$ and ${\tilde m}\in {\tilde M}$, the corresponding coefficient $c_{m,{\tilde m}}$ of $\VV$ is the complex function on $G$ defined by 
$$
c_{m,{\tilde m}}(g)=\langle {\tilde m}, \pi (g).m\rangle_\VV , \ g\in G\ .
$$
\noi Let us observe that by construction $c_{m,{\tilde m}}$ is $I$-bi-invariant. 

\begin{corollary} \label{Coefficient} For $m\in M$, ${\tilde m}\in {\tilde M}$, $w\in \WAff$, the value of $c_{m,{\tilde m}}$ on $IwI$ is $q_w^{-1}\langle {\tilde m}, r({e}_w ) .m\rangle_M$. 
\end{corollary}
\noi {\it Proof}.  The vector $m\in \VV$ corresponds to the function $v_m\in C_c (G/I ,M)$ which has support $I$ and constant value $m$ on $I$. The linear form $\tilde m$ corresponds to the function $f_{\tilde m}\in C(G/I,{\tilde M})$ whose value on $IwI$, $w\in \WAff$, is given by $q_w^{-1}\, {\tilde r}(\check{e}_w).m$.  
 Applying the pairing formula (\ref{Pairing}), for $w\in \WAff$, we have $c_{m,{\tilde m}}(w) = \langle q_w^{-1}\, {\tilde r}(\check{e}_w ).{\tilde m},m\rangle_M = q_w^{-1}\, \langle {\tilde m}, r(e_w ).m\rangle$, as required. 
\medskip

We now describe a geometric model for the space ${\mathbf P}^o ({\tilde M})$ for a given representation 
$(r,M)$ of $\HIw$.   We let $\Ch (X_E )$ denote the $G$-set of chambers of $X_E$ and $\FF (\Ch (X_E ), {\tilde M})$ the space of $\tilde M$-valued functions on $\Ch (X_E )$. We let $\FF (\Ch (X_E ),{\tilde M})_r$ denote the subspace of $\FF (\Ch (X_E ),{\tilde M})$ formed of those functions $f$ satisfying for all $C\in \Ch (X_E )$ and for all $s\in S$:
\begin{equation}\label{RelationChambre}
{\tilde r}(e_s ). f(C) = \sum_{D\sim_s C} f(D)\ .
\end{equation}
Since $G$ acts on $\Ch (X_E )$ by preserving the labelling of vertices,  $\FF (\Ch (X_E ),{\tilde M})_r$ is a $G$-subspace of  $\FF (\Ch (X_E ),{\tilde M})_r$
\begin{proposition} \label{Model} The spaces  ${\mathbf P}^o ({\tilde M})$ and $\FF (\Ch (X_E ),{\tilde M})_r$
 are isomorphic as $G$-modules. 
\end{proposition}

\noi {\it Proof}. Since $G$ acts transitively on $\Ch (X_E )$ and the global stabilizer of $C_o$ is $I$, we may identify $C (G/I ,{\tilde M})$ with $\FF (\Ch (X_E ),{\tilde M})$. Since $\{ e_s \ ; \ s\in S\}$ generates $\HIw$ as an algebra, a function $f\in C(G/I ,{\tilde M})$ is in ${\mathbf P}^o ({\tilde M})$ if and only if $f\star e_s 
={\tilde r}(e_s ).f$, for all $s\in S$. It follows that ${\mathbf P}^o ({\tilde M})$ identifies with the space of functions $f$ :  $\Ch (X_E )\lra {\tilde M}$ satisfying $f\star e_s = {\tilde r}(e_s ).f$. Now our result follows from \cite{Bo} {\S}4 Equation (6) : 
\begin{equation}\label{Actiones}
f\star e_s (C) =\sum_{D\sim_s C} f(D), \ f\in \FF (\Ch (X_E ),{\tilde M}), \ s\in S\ .
\end{equation}

\section{Multiplicities}

 We continue to abreviate $I=I_E^o$, a fixed ${\rm Gal}(E/F)$-stable Iwahori subgroup of $G$.

  We fix an irreducible  Iwahori-spherical representation $(\pi , \VV )$ of $G$. It is  of the form ${\mathbf V}(M)$, for some irreducible representation $(r,M)$ of $G$.  The aim of this section is to give a bound for the dimension of the intertwining space ${\rm Hom}_H \,  (\VV ,\CC )$, where $\CC$ is acted upon via the trivial representation of $H$.
 
We shall prove the following. 

\begin{theorem} \label{Multiplicite}  The obvious restriction map ${\rm Hom}_H (\pi ,\CC )\lra {\rm Hom}_\CC (\pi^I ,\CC )$ is injective. In particular, we have
$$
{\rm dim}_\CC\, {\rm Hom}_H \, ({\mathbf V}(M), \CC )\leqslant {\rm dim}_\CC\, M\ .
$$
\end{theorem}

 By Propositions \ref{Dual} and \ref{Model},  we have an isomorphism of vector spaces:
$$
  {\rm Hom}_H \, ({\mathbf V}(M), \CC )= ({\mathbf V}(M)^* )^H \simeq {\mathbf P}^o ({\tilde M})^H
\simeq \FF (\Ch (X_E ),{\tilde M})_r^H \ .
$$
\noi  Through this isomorphism, the restriction map ${\rm Hom}_H (\pi ,\CC )\lra {\rm Hom}_\CC (\pi^I ,\CC )$ corresponds to the map
$\FF (\Ch (X_E ),{\tilde M})_r^H \lra {\tilde M}$, $f\lra f(C_o )$. 
 Hence we are reduce to proving the following result.

\begin{proposition} \label{Injection} A map $f\in \FF (\Ch (X_E ),{\tilde M})_r^H$ is entirely determined by the value $f(C_o )$, i.e. the mapping $\FF (\Ch (X_E ),{\tilde M})_r^H \lra {\tilde M}$, $f\mapsto f(C_o )$ is injective.
\end{proposition}

 To prove the proposition we use the same tool that  allowed to establish the multiplicity $1$ property of the main result  of \cite{BC} (Theorem 2):

\begin{proposition} \label{Transitive}  (Fran\c cois Court\`es). Let $d$ be a non-negative integer and $C$ be a chamber of $X_E$ at combinatorial distance $d$ from $X_F$. Let $D^+$ be a chamber of $X_E$, adjacent to $C$ and at combinatorial distance $d+1$ from $X_F$. Consider the wall $A=D^+ \cap C$. The set ${\mathcal C}_A$ of chambers of $X_E$ containing $A$ splits into two subsets: the set ${\mathcal C}_A^+$ (resp. ${\mathcal C}_A^-$) formed of those chambers at combinatorial distance $d+1$ (resp. $d$) from $X_F$.

 Then the stabilizer $H_A$ of $A$ in $H$ acts transitively of ${\mathcal C}_A^+$ et ${\mathcal C}_A^-$. Moreover we have :
$$
(\vert {\mathcal C}_A^-\vert ,\vert {\mathcal C}_A^+\vert )=\left\{ \begin{array}{l}
(1,q_o^2 )\\
{\rm or}\\
(q_o +1 ,q_o^2 -q_o )
\end{array}\right.
$$
\end{proposition}

\noi {\it Proof}. This follows from  \cite{BC} Prop. A1 and the proof of \cite{Cou} Theorem 6.1.
\medskip

\noi {\it Proof of Proposition \ref{Injection}}.   Let $f\in \FF (\Ch (X_E ),{\tilde M})_r^H$. For a non-negative integer $d$, we let $\Ch (X_E )_d$ denote the set of chambers in $X_E$ at combinatorial distance from $X_F$ less than or equal to  $d$. In particular $\Ch (X_E )_0 =\Ch (X_F )$,  the set of chambers in $X_F$.  We prove by induction on $d$,  that  for all $d\geqslant 0$, the restriction of $f$ to $\Ch (X_E )_d$ depends only on $f(C_o )$. 

 Our assertion holds for $d=0$ since $H$ acts transitively on $\Ch (X_F )$. Assume it holds for an integer $d\geqslant
0$.  Let $D$ be a chamber a combinatorial distance $d+1$ from $X_F$. By definition there exists a chamber $C$ adjacent to $D$ and at combinatorial distance $d$ from $X_F$. Use the notation  $A=D\cap C$, ${\mathcal C}_A^-$ and ${\mathcal C}_A^+$  as in Proposition \ref{Transitive}. Let $s$ be the type of $A$.  By Relation \ref{RelationChambre} we have

\begin{equation}
{\tilde r}(e_s ).f(C)=\sum_{E\in {\mathcal C}_A^- \backslash \{ C\}} f(E) +\sum_{E\in {\mathcal C}_A^+} f(E)\\
= (\vert {\mathcal C}_A^- \vert -1)f(C) +\vert {\mathcal C}_A^+\vert f(D)
\end{equation}

\noi where we used the fact that $H$ acts transitively on ${\mathcal C}_A^+$ and ${\mathcal C}_A^-$, and that $f$ is $H$-invariant. Hence we have
\begin{equation}
f(D) =\frac{1}{\vert {\mathcal C}_A^+ \vert} {\tilde r}\left( e_s  -(\vert {\mathcal C}_A^-\vert -1)\, e_1 \right) . f(C)
\end{equation}

\noi so that $f(D)$ depend only on $f(C_o )$, as required.

\medskip

Theorem \ref{Multiplicite} has the following nice consequence. 

\begin{theorem} \label{Mult1} Let $(\pi, \VV )$ be an irreducible Iwahori-spherical representation satisfying ${\rm dim}\, \VV^I =1$. Then $\pi$ has the multiplicity $1$  property:
$$
{\rm dim}_\CC \,{\rm Hom}_H\, (\VV ,\CC )\leqslant 1\ .
$$
\end{theorem}

\section{Test vectors}

 Let us fix an irreducible square integrable (or discrete series) representation $(\pi ,\VV )$ of $G$. By definition, for $v\in \VV$ and ${\tilde v}\in {\tilde \VV}$, the coefficient $c_{v,{\tilde v}}$ : $G\ni g\mapsto \langle {\tilde v}, \pi (g).v\rangle_\VV$ lies in $L^2 (G)$.  Fix a Haar measur $\nu$ on $H$. 

\begin{proposition}\label{integrable} For all $v\in \VV$, ${\tilde v}\in {\tilde \VV}$, the integral
$$
\int_{H} c_{v,{\tilde v}}(h)\, d\nu (h)
$$
\noi is absolutely convergent. 
\end{proposition}

\noi {\it Proof}. Indeed, by \cite{GO} Corollary 1.2, a Galois symmetric space as $G/H$ is strongly discrete: the restriction to $H$ of a coefficient of an irreducible discrete series representation of $G$ lies in $L^1 (H)$. 
\medskip
 
 Following Zhang \cite{Zh}, for $v\in \VV$ and ${\tilde v}\in {\tilde \VV}$, we set
$$
{\mathcal L} (v,{\tilde v}) = \int_{H} c_{v,{\tilde v}}(h)\, d\nu (h)\ .
$$
\noi Then $\mathcal L$ : $\VV\times {\tilde \VV}\lra \CC$ is a bi-$H$-invariant bilinear form. For ${\tilde v}\in {\tilde \VV}$, we set
$$
{\mathcal L}_{\tilde v}\ : \ \VV\lra \CC , \ v\mapsto {\mathcal L}(v,{\tilde v})\ .
$$
\noi By construction ${\mathcal L}_{\tilde v}\in {\rm Hom}_H\, (\VV ,\CC )$, for all ${\tilde v}\in {\tilde \VV}$ and we set ${\mathcal H}(\pi ) = \{ {\mathcal L}_{\tilde v} \ ; \ {\tilde v}\in {\tilde \VV}\}\subset {\rm Hom}_H \, (\VV ,\CC )$. 

 By \cite{Zh} Proposition 3.2, $G/H$ is very strongly discrete in the sense of {\it loc. cit.} Definition 1.3; this means that the linear form
$$
{\mathcal C}(G)\lra \CC , \ f\mapsto \int_{H} f(h)\, d\nu (h)
$$

\noi is well-defined and continuous, where ${\mathcal C}(G)$ denotes the  Harish-Chandra's Schwartz space of $G$ (cf. {\it loc. cit.} {\S}2).  As a consequence, by {\it loc. cit.} Theorem 1.4, we have the following result. 

\begin{theorem} \label{Egalite} (Zhang) We have the equality ${\mathcal H}(\pi )={\rm Hom}_H \, (\VV ,\CC )$, i.e. any $H$-invariant linear form on $\VV$ is obtained by integration on $H$ against a smooth linear form on $\VV$. 
\end{theorem}

 The following result relates the distinction of $\pi$ to the integral of a bi-$I$-invariant coefficient. 

\begin{theorem} \label{Bitest} The discrete series representation $(\pi ,\VV )$ is $H$-distinguished if, and only if, there exist $v\in \VV^I$ and ${\tilde v}\in {\tilde\VV}^I$ such that $\ds \int_{H}c_{v,{\tilde v}}(h)\, d\nu (h)\not= 0$. 
\end{theorem}

 {\it Proof}.  One implication in Theorem \ref{Bitest} being obvious, let us assume that $(\pi, \VV )$ is $H$-distinguished and let us fix a non-zero $\Lambda\in {\rm Hom}_H (\VV ,\CC )$.  By Theorem \ref{Egalite} there exists ${\tilde v}_1\in {\tilde \VV}$ such that $\Lambda ={\mathcal L}_{{\tilde v}_1}$. By Theorem \ref{Multiplicite} pick $v\in \VV^I$ such that $\Lambda (v)={\mathcal L}(v,{\tilde v}_1 )\not= 0$. The map ${\tilde \Lambda}$ : ${\tilde v}\mapsto {\mathcal L}(v,{\tilde v})$ is a non-zero element of ${\rm Hom}_H \, ({\tilde\VV},\CC )$. Applying  Theorem \ref{Multiplicite} again to $({\tilde \pi},{\tilde \VV})$, there exists ${\tilde v}\in {\tilde \VV}^I$ such that ${\mathcal L}(v,{\tilde v})\not= 0$ and we are done.

\section{Generalized Poincar\'e series}

Recall that we abbreviate  $W=W_F^{\rm Aff}=W_E^{\rm Aff}$. 

 The Poincar\'e series of the Coxeter system $(W,S)$ is the formal power series $W(t)=\ds\sum_{w\in W} t^{l(w)}\in \ZZ [[t]]$. This is the generating  function of the length of elements of $W$. It is known (\cite{Bott},\cite{St}) to be a rational function of $t$. We shall need a generalization of the Poincar\'e series due to Macdonald and Gyoja. 
\medskip

 We let $S_1$, ..., $S_m$ be the non-empty intersections of $S$ with the conjugacy classes of $W$. By \cite{Bki} Proposition 5, page 16, for $w\in W$ the number of elements of $S_i$ occuring in a reduced expression of $w$ depends only $w$, we denote it by $l_i (w)$. We have $l(w) = l_1 (w)+\cdots +l_m (w)$. Let $t_1$, ..., $t_m$ be indeterminates. For $w\in W$, we write ${\mathbf t}=(t_1 ,...,t_m )$ and  ${\mathbf t}^{l(w)} := t_1^{l_1 (w)} \cdots t_m^{l_m (w)}$. 

  By \cite{Bo} {\S}3.3, the $S_i$, $i=1,...,m$ are the connected components of the subgraph of the Coxeter graph ${\rm Cox}(W,S)$ obtained by erasing the multiple edges. The number $m$ is so given as follows:

\smallskip

 $m=1$ if $G$  is of type $A_n$ ($n\geqslant 2$), $D_n$ ($n\geqslant 3$) and $E_i$ ($i=6,7,8$), 

 $m=2$ if $G$ is of type $A_1$, $B_n$ ($n\geqslant 2$), $G_2$ and $F_4$,

$m=3$ if $G$ is of type $C_n$ ($n\geqslant 2$). 
\smallskip

\noi  We choose the indexing in $(t_i )_{i=1,..,m}$ as in \cite{Gy}, pages 173, 174.

\medskip
 
 Following Gyoja \cite{Gy}, for a representation $(r,M)$ of $\HIw$, define an element of ${\rm End}_\CC (M)\otimes_\CC \CC [[t_1 ,...,t_m ]]$ by 
\begin{equation} \label{Gyoja}
L({\mathbf t},r) = \sum_{w\in W} r(e_w )\, {\mathbf t}^{l(w)}\ .
\end{equation}

 \begin{theorem} (\cite{Gy}   In any basis of $M$, the coefficients of $L({\mathbf t},r)$ are rational functions of $t_1$, ..., $t_m$.
\end{theorem} 

When $(r,M)$ is the trivial representation of $\HIw$, we set $W ({\mathbf t})=L({\mathbf t},r)$, that is 
$W({\mathbf t}) = \ds \sum_{w\in W} {\mathbf t}^{l(w)}$. This rational function has been computed by Macdonald (\cite{Mc} Theorem 3.3). When $m=1$, then $W({\mathbf t})= W(t)$ is the ordinary Poincar\'e series and Bott's fomula gives
$$
W(t) =   \prod_{i=1,..,l} \frac{1-t_i^{m_i +1}}{(1-t_i)(1-t_i^{m_i})}
$$
\noi where $m_1$, ..., $m_l$ are the exponents of the spherical Coxeter group attached to $W$ (cf. \cite{Bki} Chap. V, {\S}6, D\'efinition 2). 
\smallskip

 If $m=2$, the formulas are the following.  
\smallskip

  Type $A_1$, $W(t_1 ,t_2 )=\ds \frac{(1+t_1 )(1+t_2 )}{1-t_1 t_2 }$.
 
  Type $B_n$ ($n\geqslant 2$), $\ds W(t_1 ,t_2 )=\frac{1-t_1^n }{(1-t_1 )^n}\, \prod_{i=1}^{n-1} (1-t_1^{2i})\, \prod_{i=0}^{n-1} \frac{ 1+t_1^i t_2}{1-t_1^{n-1+i}t_2}$.

  Type $G_2$, $\ds W (t_1 ,t_2 )=\frac{(1+t_1 )(1+t_1 +t_1^2 )(1+t_2 )(1+t_1 t_2 +t_1^2 t_2^2 )}{(1-t_1^2 t_2 )(1-t_1^3 t_2^2 )}$.
 
 Type $F_4$, $W(t_1 ,t_2 )=$
 
$$ \left( \prod_{i=1}^3 \frac{(1-t_1^{i+1})(1+t_1^i t_2 )(1-t_2^i )}{(1-t_1 )(1-t_2 )}\right) \, \frac{(1+t_1 t_2^2 )(1+t_1^2 t_2^2 )(1+t_1^3 t_2^3 )}{(1-t_1^3 t_2^2 )(1-t_1^4 t_2^3 )(1-t_1^5 t_2^3 )(1-t_1^6 t_2^5 )} \ .
$$
 \smallskip

 Finally for $m=3$ we have the following formula.
\smallskip

 Type $C_n$ ($n\geqslant 2$), $\ds W(t_1 ,t_2 ,t_3 ) = \prod_{i=0}^{n-1}
\frac{(1-t_1^{i+1})(1+t_1^i t_2 )(1+t_1^i t_3 )}{(1-t_1 )(1-t_1^{n-1+i} t_2 t_3 )}$.

\section{A numerical criterion with application to degree $1$ characters}

 Let $(\pi ,\VV )$ be an irreducible discrete series representation of $G$ and let us assume that it is Iwahori-spherical: $\VV^I\not= 0$. Let $(r,M)$ denote the representation of $\HIw$ in $M=\VV^I$. The following result is a numerical criterion to decide whether or not $(\pi ,\VV )$ is $H$-distinguished. 

\begin{theorem}\label{Criterion} (a) The representation $\pi$ is $H$-distinguished if, and only if, the endomorphism $L({\mathbf t}_o ,r)\in {\rm End}_\CC (M)$ is non-zero, where ${\mathbf t}_o =\ds (\frac{1}{q_o},...,\frac{1}{q_o})$.

(b) We have ${\rm dim}\, {\rm Hom}_H \, (\pi ,\CC )\geqslant {\rm rank}\, L({\mathbf t}_o ,r)$.
\end{theorem}

\noi {\it Proof}. (a) We abbreviate $J=I_F^o$. Recall that $JwJ\subset IwI$, for all $w\in W$. Let $\nu$ the Haar measure on $H$ normalized by $\nu (J)=1$. Let $m\in M$ and ${\tilde m}\in {\tilde M}$. We first calculate the integral $\ds\int_H c_{m,{\tilde m}} (h)\, d\nu (h)$. By the Bruhat-Tits decomposition $\ds H=\sqcup_{w\in W} JwJ$, we have:
\begin{align}
 \int_{H} c_{m,{\tilde m}} (h)\, d\nu (h) & = \sum_{w\in W} \int_{JwJ} c_{m,{\tilde m}}(h)\, d\nu (h)\\
 \label{*} & = \frac{\nu (JwJ)}{q^{l(w)}} \sum_{w\in W} \langle {\tilde m}, r(e_w ).m\rangle_M\\
& = \frac{q_o^{l(w)}}{q^{l(w)}} \sum_{w\in W}\langle {\tilde m}, r(e_w ).m\rangle_M\\
& = \langle {\tilde m} ,\sum_{w\in W}r(e_w )\left( \frac{1}{q_o}\right)^{l(w)} .m\rangle_M\\
 & =\langle {\tilde m}, L({\mathbf t}_o ,r).m\rangle_M
\end{align}
\noi where in (\ref{*}) we used Corollary \ref{Coefficient} and where we set ${\mathbf t}_o = \ds (\frac{1}{q_o},...,\frac{1}{q_o})$.

By Proposition \ref{Bitest}, $\pi$ is $H$-distinguished if and only if there exist $m\in M$ and ${\tilde m}\in {\tilde M}$ such that $\ds\int_{H} c_{m,{\tilde m}}(h)\, d\nu (h)\not= 0$. Thanks to the last calculation this is indeed equivalent to $L({\mathbf t}_o ,r)\not= 0$, as required. 
\smallskip

(b) It follows from (a), using the notation of {\S}5, that for $m\in M$ and ${\tilde m}\in {\tilde M}$, we have ${\mathcal L}_{\tilde m}(m)=\langle {\tilde m}, L({\mathbf t}_o ,r)\rangle_M$.  Let us put $s={\rm rank}\, L({\mathbf t}_o ,r)$. Then there exist ${\tilde m}_1$, ..., ${\tilde m}_s\in{\tilde M}$ such that the $({\mathcal L}_{{\tilde m}_i})_{\mid M}$, $i=1,...,s$, are linearly independent. It follows that the ${\mathcal M}_{{\tilde m}_i}\in {\rm Hom}_H \, (\pi ,\CC )$, $i=1,...,s$, are linearly independent, and we are done.

\medskip

 We now assume that the discrete series representation $(\pi ,\VV )$ satisfies ${\rm dim}\, M = {\rm dim}\, \VV^I =1$, so that we may view $r$ as an algebra homomorphism $\HIw \lra \CC$. Such discrete series where classified in \cite{Bo} {\S}5 by A. Borel. We recall this classification. 

 The character $r$ is constant on each $\{ e_s \ ; \ s\in S_i\}$, $i=1,...,m$ and because of the quadratic relations we have $r(e_s )\in \{ -1 ,q\}$, $s\in S$. Following Borel {\it loc. cit},  for $i=1,...,m$ we set
$$
\epsilon_i =\left\{ \begin{array}{rl}
1 & {\rm if}\ r(e_s )=q\\
-1 & {\rm if}\ r(e_s )=-1
\end{array}\right.
$$ 
\noi so that the representation $\pi$ is entirely characterized by the $m$-uple $(\epsilon_1 ,...,\epsilon_m )$; we write $\pi = \pi (\epsilon_1 ,...,\epsilon_m )$. Recall that $\St_G = \pi (-1,..,-1 )$ is the Steinberg 
representation of $G$. 

 \begin{proposition} \label{List} (Borel \cite{Bo} {\S}5.8) The discrete series representations of $G$ of the form $\pi (\epsilon_1 ,...,\epsilon_m )$ are the Steinberg representation and

$\pi (-1,1)$ for types $B_n$ ($n\geqslant 3$), $F_4$, $G_2$,

$\pi (-1,-1,1)$, $\pi (-1,1,-1)$ for types $C_2$, $C_3$,

$\pi (-1,-1,1)$, $\pi (-1,1,-1)$, $\pi (-1,1,1)$ for type $C_n$ ($n\geqslant 4$). 
\end{proposition}

 \begin{lemma} \label{Value} Assume that $\pi (\epsilon_1 ,...,\epsilon_m )$ is a discrete series representations of $G$. We have $L({\mathbf t}_o ,r)= W(\epsilon_1 q_o^{\epsilon_i},...,\epsilon_m q_o^{\epsilon_m})$. In particular $\pi (\epsilon_1 ,...,\epsilon_m )$ is $H$-distinguished if, and only if, $W(\epsilon_1 q_o^{\epsilon_i},...,\epsilon_m q_o^{\epsilon_m}) \not= 0$. 
\end{lemma}
\noi {\it Proof}. Observe that for $i=1,...,m$ and $s\in S_i$, we have $r(e_s )= \epsilon_i q_{o}^{\epsilon_i +1}$, so that for $w\in W$ we have $\ds r(e_w ) = \prod_{i=1}^m (\epsilon_i q_o^{\epsilon_i +1})^{l_i (w)}$. Hence we have
\begin{align*}
L({\mathbf t}_o ,r) &= \sum_{w\in W} r(e_w ) (\frac{1}{q_o})^{l(w)}\\
 &= \sum_{w\in W}\prod_{i=1}^m (\epsilon_i \frac{q_o^{\epsilon_i +1}}{q_o})^{l_i (w)}\\
 & = W( \epsilon_i q_o^{\epsilon_i}, ..., \epsilon_i q_o^{\epsilon_i})
\end{align*}
\noi as required. 

We may now state the main result of this article. 

\begin{theorem} \label{Main} The discrete series representation $\pi (\epsilon_1 ,...,\epsilon_m )$ is $H$-distinguished if and only if it is a Steinberg representation, or if $H$ is of type $G_2$ and $(\epsilon_1 ,\epsilon_2 )=(-1 ,1)$. 
 When it is distinguished it has the multiplicity $1$ property: $\ds {\rm dim}\, {\rm Hom}_H\, (\pi (\epsilon_1 ,...,\epsilon_m ),\CC ) = 1$.
\end{theorem}

\noi {\it Proof}. The multiplicity $1$ property follows from Theorem \ref{Multiplicite}.  The first assertion is proved by a case-by-case calculation based on Gyoja's formulas given in {\S}5. Note that in the case of the Steinberg representation, the distinction property was originally D. Prasad's conjecture (\cite{Pr2} Conjecture 3, p. 77) and was proved in \cite{BC} (cf. \cite{Be} for a proof of Prasad's conjecture in the general case).

Paul Broussous
\smallskip

paul.broussous{@}math.univ-poitiers.fr
\medskip

 Laboratoire de Math\'ematiques et Applications, UMR 7348 du CNRS

Site du Futuroscope -- T\'el\'eport 2
 
11,  Boulevard Marie et Pierre Curie

B\^atiment H3 --  TSA 61125

86073 POITIERS CEDEX 

 France

\end{document}